\documentclass{article}[12pt]
\usepackage{color}
\usepackage[normalem]{ulem}
\usepackage{times}
\usepackage{fullpage}
\usepackage{amsmath}
\usepackage{amssymb}
\usepackage{verbatim}
\usepackage{color}
\usepackage{multirow}
\usepackage{paralist}
\def \R {\mathbb R}
\def \Z {\mathbb Z}
\def \N {\mathbb N}

\def \eps {\varepsilon}

\def \NE {{\sf NE}}
\def \TC {{\sf TC}}
\def \P {{\sf P}}
\def \AC {{\sf AC}}

\def \AS {\mathbf{AS}}
\def \D {\mathbf{D}}
\def \bit {\lbrace -1, 1 \rbrace}
\def \booleanfn {\hypercube \rightarrow \bit}
\def \hypercube {\bit^n}
\def \sgn {\text{sgn}}

\def \Lemming {Lemma}

\def \cm {HSF}
\def \cms {HSFs}
\def \acm {an HSF}
\def \tick {$\color{cyan}\checkmark$}
\def \cross {$\color{red}\times$}
\newenvironment{proof}{{\bf Proof.  }}{\hfill$\Box$}

\newtheorem{theorem}{Theorem}[section]
\newtheorem{definition}{Definition}[section]
\newtheorem{corollary}{Corollary}[section]
\newtheorem{lemma}{\Lemming}[section]

\newtheorem{conjecture}{Conjecture}[section]
\title{The Gotsman-Linial Conjecture is False}
\author{Brynmor Chapman\thanks{Supported by an NSF Graduate Research Fellowship}\\MIT}
\date{1 April 2017}

\begin{document}

\maketitle

\begin{abstract}
In 1991, Craig Gotsman and Nathan Linial conjectured that for all $n$ and $d$, the average sensitivity of a degree-$d$ polynomial threshold function on $n$ variables is maximized by the degree-$d$ symmetric polynomial which computes the parity function on the $d$ layers of the hypercube with Hamming weight closest to $n/2$.
We refute the conjecture for almost all $d$ and for almost all $n$, and we confirm the conjecture in many of the remaining cases.
\end{abstract}

\section{Introduction}

We say that a boolean function $f : \booleanfn$ is a \emph{Polynomial Threshold Function of degree $d$} if it can be expressed as the sign of a polynomial $p \in \R [ x_1, \ldots, x_n ]$ of degree at most $d$ evaluated on the boolean hypercube.
For brevity, we will use the term $(n,d)$-PTF (or simply PTF, when $n$ and $d$ are either implicit or irrelevant) to refer to a polynomial threshold function of degree $d$ on $n$ variables.
We say that the coefficients of $p$ are the \emph{realizing weights} of $f$.
Note that these realizing weights are not unique, as any sufficiently small perturbation of $p$ will not affect its sign on the discrete set $\hypercube$.
This definition alone is not terribly exciting without restrictions on $d$, as every boolean function on $n$ variables can be written as the sign of (and in fact can be written exactly as) a multilinear polynomial of degree $n$.
We are interested particularly in the case where $d$ is small.

In an influential paper, Craig Gotsman and Nathan Linial \cite{GL94} applied Fourier analytic techniques to the study of PTFs.
They were mainly interested in connecting different measures of the complexity of boolean functions, and of low-degree PTFs in particular.
One such measure was the Average Sensitivity of a boolean function, defined in Fourier analytic terms.
For simplicity, in this paper we use the following (equivalent) combinatorial definition:
\begin{definition}
For a function $f: \booleanfn$, we define its \emph{Dichromatic Count} $\D[f]$ to be the number of (unordered) pairs of Hamming neighbors $\{x,y\}$ such that $f(x) \ne f(y)$.
\end{definition}
We say that such a pair of Hamming neighbors is a \emph{dichromatic edge} of $f$.
\begin{definition}
The \emph{Average Sensitivity} of a boolean function $f$ is $\AS[f] := 2^{1-n} \D[f]$.
\end{definition}

Among other things, Gotsman and Linial proved a tight upper bound on the average sensitivity of $(n,1)$-PTFs, achieved by the MAJORITY function on $n$ variables.
They conjectured that this bound generalizes to higher degree PTFs, in that the $(n,d)$-PTF of maximal average sensitivity is the obvious symmetric candidate, which alternates signs on the $d+1$ values of $\displaystyle \sum_{i \in [n]} x_i$ closest to $0$.
\begin{conjecture}[Gotsman-Linial]\label{GL}
Let $p^*_{n,d}$ be the monic univariate polynomial of degree $d$ with (non-repeated) roots at the $d$ integers closest to $0$ of opposite parity from $n$.
Let $\displaystyle f^* ( x_1, \ldots, x_n ) = \sgn \left( p^*_{n,d} \left( \sum_{i \in [n]} x_i \right) \right)$.
Then for every $(n,d)$-PTF $f$, $\AS[f] \le \AS[f^*_{n,d}]$.
\end{conjecture}
This conjecture was listed as a prominent open problem in \cite{OD14} and \cite{alphabetsoup}.
If true, it would have many applications in complexity and learning (see for example \cite{HKM, GS, K12, KW, CSS}), although most of the applications would already be implied by an asymptotic version of the conjecture, stated below.
Gotsman and Linial proved their conjecture for the case where $d=1$, and it is also known to be true in the case where $d=0$.
However, it was left open whether the conjecture holds for any $d \ge 2$.
Two weaker versions of this conjecture have since been formulated and studied.
\begin{conjecture}[Gotsman-Linial - Asymptotic]\label{GLweak}
Let $f : \booleanfn$ be an $(n,d)$-PTF.
Then the average sensitivity $\AS[f] \in O( d \sqrt{n} )$.
\end{conjecture}
\begin{conjecture}[Gotsman-Linial - Weak]\label{GLsuperweak}
Let $f : \booleanfn$ be an $(n,d)$-PTF.
Then the average sensitivity $\AS[f] \in O( \sqrt{n} \log^{g(d)} n )$ for some function $g$ depending only on $d$.
\end{conjecture}
Conjecture \ref{GLsuperweak} was resolved by Daniel Kane \cite{K13}.

\subsection{Result}

In this paper, we resolve the Gotsman-Linial Conjecture (Conjecture \ref{GL}) for all pairs $(n, d)$ except the case when $n>7$ is even and $d=2$.
The main result of this paper is the following.

\begin{theorem}\label{mainrefute}
For all pairs of natural numbers $(n,d)$ satisfying one of the following criteria, there exists an $(n,d)$-PTF $f_{n,d}$ witnessing a counterexample to the Gotsman-Linial Conjecture (Conjecture \ref{GL}):
\begin{compactitem}
\item $n \ge 5$ is odd, and $d=2$.
\item $n \ge 7$, and $3 \le d \le n-3$.
\end{compactitem}
Moreover, $\AS[f_{n,d}] \in (1 + \Omega(n^{-1} e^{-d^2/n})) \AS[f^*_{n,d}]$.
\end{theorem}

In addition, the conjecture holds in many of the remaining cases.

\begin{theorem}\label{mainconfirm}
For all pairs of natural numbers $(n,d)$ satisfying one of the following criteria, $f^*_{n,d}$ has the greatest average sensitivity among $(n,d)$-PTFs.
\begin{compactitem}
\item $d \le 1$.
\item $d \ge n-2$.
\item $n=6$.
\end{compactitem}
\end{theorem}

Our results (and the remaining open cases) are summarized in Figure \ref{resultstable}.
Although we refute the Gotsman-Linial Conjecture for most cases that are of interest for applications, the asymptotic conjecture (Conjecture \ref{GLweak}), which would suffice for most known applications, remains open.

\begin{figure}[h]
\begin{center}
\begin{tabular}{cr|cccccccccccccccc}
& \multicolumn{1}{r}{} & \multicolumn{16}{c}{$n$} \\
&& $1$ & $2$ & $3$ & $4$ & $5$ & $6$ & $7$ & $8$ & $9$ & $10$ & $11$ & $12$ & $\cdots$ & $2k$ & $2k+1$ & $\cdots$ \\
\cline{2-18}
\multirow{9}{*}{$d$} & $0$ & \tick & \tick & \tick & \tick & \tick & \tick & \tick & \tick & \tick & \tick & \tick & \tick & $\color{cyan} \cdots$ & \tick & \tick & $\color{cyan} \cdots$ \\
& $1$ & \tick & \tick & \tick & \tick & \tick & \tick & \tick & \tick & \tick & \tick & \tick & \tick & $\color{cyan} \cdots$ & \tick & \tick & $\color{cyan} \cdots$ \\
& $2$ && \tick & \tick & \tick & \cross & \tick & \cross & ? & \cross & ? & \cross & ? & $\cdots$ & ? & \cross & $\cdots$ \\
& $3$ &&& \tick & \tick & \tick & \tick & \cross & \cross & \cross & \cross & \cross & \cross & $\color{red} \cdots$ & \cross & \cross & $\color{red} \cdots$ \\
& $4$ &&&& \tick & \tick & \tick & \cross & \cross & \cross & \cross & \cross & \cross & $\color{red} \cdots$ & \cross & \cross & $\color{red} \cdots$ \\
& $5$ &&&&& \tick & \tick & \tick & \cross & \cross & \cross & \cross & \cross & $\color{red} \cdots$ & \cross & \cross & $\color{red} \cdots$ \\
& $6$ &&&&&& \tick & \tick & \tick & \cross & \cross & \cross & \cross & $\color{red} \cdots$ & \cross & \cross & $\color{red} \cdots$ \\
& $7$ &&&&&&& \tick & \tick & \tick & \cross & \cross & \cross & $\color{red} \cdots$ & \cross & \cross & $\color{red} \cdots$ \\
& $\vdots$ &&&&&&&& $\color{cyan} \ddots$ & $\color{cyan} \ddots$ & $\color{cyan} \ddots$ & $\color{red} \ddots$ & $\color{red} \ddots$ & $\color{red} \ddots$ & $\color{red} \ddots$ & $\color{red} \ddots$ & $\color{red} \ddots$
\end{tabular}
\end{center}
\caption{Results are summarized in the above table.
A cyan tick mark indicates a case in which the conjecture holds (for all $(n,d)$-PTFs $f$, $\AS[f] \le \AS[f^*_{n,d}]$).
A red cross indicates a refutation (there exists an $(n,d)$-PTF $f$ such that $\AS[f] \in (1 + \Omega(n^{-1} e^{-d^2/n})) \AS[f^*_{n,d}]$).
A black question mark indicates an open case.
Note: the cases $(n,d) = (6,2)$ and $(n,d) = (6,3)$ were verified with the help of a computer search and a linear program solver (see Appendix~\ref{appendix-computer}).}
\label{resultstable}
\end{figure}

The remainder of this paper is structured as follows.
We first present some high level intuition relating to the Gotsman-Linial Conjecture.
Section 2 contains background information.
Section 3 contains constructions of the refutations indicated in Figure \ref{resultstable}.
Section 4 concludes the paper and presents a revised conjecture.

\subsection{Intuition}\label{intuition}
We start with some very high level intuition as to why the Gotsman-Linial Conjecture might be (approximately) true.
The conjecture holds in the case of symmetric PTFs (boolean functions which can be expressed as the sign of a univariate polynomial in the sum of the input bits).
This follows from the Fundamental Theorem of Algebra and a simple counting argument.
In the more general case, we might expect that a degree-$d$ PTF can be expressed (at least approximately) in terms of $d$ unate functions.
This generalizes the observation that every linear threshold function is unate.
For a sufficiently close approximation, this would prove the Asymptotic Gotsman-Linial Conjecture.
Intuition may also be drawn from Kane's proof of Conjecture \ref{GLsuperweak}.
If inputs are chosen from a Gaussian distribution instead of a Bernoulli distribution, a polynomial $p$ is expected to be too large in magnitude for a small change in its input to change its sign.
Under certain conditions, a similar result can be extended to polynomial threshold functions on the boolean hypercube.

As for why the Gotsman-Linial Conjecture is not (exactly) true, we observe that the PTF of conjectured maximal average sensitivity is the product of $d$ linear threshold functions, with parallel separating hyperplanes between two of the middle $d+1$ layers (sets of vertices of equal Hamming weight) in the hypercube.
For some $d$, one might expect to be able to find a PTF of greater average sensitivity approximated by turning one of these separating hyperplanes `sideways', i.e. replacing a hyperplane that cuts the fewest edges with a hyperplane orthogonal to the rest.
Intuitively, this would require that $d$ be sufficiently large that some of the hyperplanes cut many more edges than others, but also sufficiently small that not too many edges are cut by two hyperplanes.
As it turns out, this intuition can be formalized for many $n$ and $d$, refuting the Gotsman-Linial Conjecture.

\section{Preliminaries}

\subsection{Background}
Low-degree PTFs, in particular linear threshold functions (degree-$1$ PTFs) with integral and polynomially bounded realizing weights, are of interest in the study of complexity classes such as $\TC$ (i.e. circuits composed of AND, OR, NOT, and MAJORITY gates of unbounded fan-in) and of neural networks.
More generally, we say that a circuit (with unbounded fan-in) is a \emph{degree-$d$ polynomial threshold circuit} if each of its constituent gates computes a degree-$d$ PTF of its inputs.
Note that since AND, OR, MAJORITY, and NOT are all linear threshold functions, $\AC$ and $\TC$ circuits are degree-$1$ polynomial threshold circuits.
Despite much research, the power of polynomial threshold circuits is poorly understood.
For instance, it is currently an open question, and a rather embarrassing one at that, whether $\NE$ (the class of functions computable in nondeterministic $2^{O(n)}$ time) is contained in $\TC^0_3$ (the class of functions computable by families of depth-$3$, polynomial size linear threshold circuits with polynomially bounded realizing weights).
Recent work by Daniel Kane and Ryan Williams \cite{KW} gave a partial answer to this question.
They studied the sensitivity of PTFs to random restrictions, proving (among other things) that $\NE$ (and in fact, $\P$-uniform $\TC^0$) does not have depth-3 $\TC$ circuits of $n^{1.499}$ gates or $n^{2.499}$ wires.

\subsection{Progress}

Conjecture \ref{GL} is trivially true in the cases $d=0$ and $d=n$ (the only $(n,0)$-PTFs are the constant functions, and $f^*_{n,n}$ is the parity function, which has the maximum possible average sensitivity).
Gotsman and Linial originally noted that Conjecture \ref{GL} had already been proven in the case where $d=1$ by Patrick O'Neil in 1971 \cite{ON}.

\begin{theorem}[O'Neil]\label{O'Neil}
The maximal number $k$ of edges of $H := \hypercube$ which may be cut by a hyperplane $P$ is given by $\displaystyle k = \left( n - \left\lfloor \frac{1}{2} n \right\rfloor \right) {n \choose \frac{1}{2} n}$.
\end{theorem}

Very little additional progress was made towards resolving the above conjectures until recently.
The first non-trivial bounds on the average sensitivity of PTFs of arbitrary degree were found independently by two groups \cite{HKM, DRST} and published jointly \cite{souplite}.
Daniel Kane in 2012 obtained the first bound which was truly sublinear in $n$ \cite{K12}, and in 2013, he proved the weak version of the Gotsman-Linial Conjecture (Conjecture \ref{GLsuperweak}) \cite{K13}.

\section{Resolution of Gotsman-Linial Conjecture}

For simplicity, we start by introducing some notation.
\begin{definition}
Let $f, g : \booleanfn$.
We say $f \sim g$, iff there exist $\sigma \in S_n$ and $\alpha \in \hypercube$ such that the function $x \mapsto f \left( x_1, \ldots, x_n \right) g \left( \alpha_1 x_{1 \sigma}, \ldots, \alpha_n x_{n \sigma} \right)$ is a constant.
\end{definition}
Note that $\sim$ defines an equivalence relation on boolean functions.
Two functions are equivalent iff one can be turned into the other through a combination of permuting the inputs and negating the inputs/output.
\begin{definition}
An $(n,d)$-Hypersensitive Function, or {\em $(n,d)$-\cm{}} is an $(n,d)$-PTF $f$ such that $\D[f] > \D[f^*_{n,d}]$.
\end{definition}
More generally, we say that a PTF $f$ is \acm{} if $n$ and $d$ are either implicit or irrelevant.
We may now restate the original Gotsman-Linial Conjecture (Conjecture \ref{GL}) as follows:

\begin{conjecture}\label{GLcm}
For all $n, d \in \N$, $(n,d)$-\cms{} do not exist.
\end{conjecture}

We first prove some simple cases of Conjecture \ref{GLcm}.
The following corollary of O'Neil's theorem (Theorem \ref{O'Neil}) uses our notation.

\begin{corollary}\label{GLd=1}
For every $n \in \N$, $(n,1)$-\cms{} do not exist.
\end{corollary}
\begin{proof}
Every $(n,1)$-PTF $f$ is defined by a separating hyperplane $P$ which cuts all of the dichromatic edges of $f$.
From O'Neil, $\displaystyle \D[f] \le ( n - \lfloor n/2 \rfloor ) {n \choose n/2} = \D[f^*_{n,1}]$, so $f$ is not \acm{}.
\end{proof}

The case $d=n-1$ is a simple consequence of a result first proven in 1968 by Marvin Minsky and Seymour Papert \cite{MP} and since re-proven several times.
We present here a variation on the proof by Aspnes et al. \cite{ABFR}.
\begin{theorem}[Minsky-Papert]
Any PTF which computes parity on $n$ variables must have degree at least $n$.
\end{theorem}
\begin{proof}
Let $p \in \R[x_1, \ldots, x_n]$ be a multilinear polynomial of degree $n-1$ which is never zero on $\hypercube$.
The set of monomials of degree at most $n$ is an orthogonal basis for the vector space of degree-$n$ multilinear polynomials on the boolean hypercube.
Hence $p$ is orthogonal to the parity function $\phi_n$, i.e. $\displaystyle \langle p, \phi_n \rangle = \sum_{x \in \hypercube} p(x) \phi_n(x) = 0$.
By assumption, every term in the sum on the RHS is non-zero, so at least one of them is negative, i.e. $\sgn \circ p \ne \phi_n$.
\end{proof}

\begin{corollary}\label{GLd=n-1}
For every $n \in \N$, $(n,n-1)$-\cms{} do not exist.
\end{corollary}
\begin{proof}
Let $f$ be an $(n,n-1)$-PTF.
Then $f \ne \phi_n$, and $f \ne -\phi_n$.
Let $X = \lbrace x : f(x) = \phi_n(x) \rbrace$ and $Y = \lbrace y : f(y) \ne \phi_n(y) \rbrace$.
Take $x \in X$ and $y \in Y$.
There are $n$ edge-disjoint paths between $x$ and $y$ in the boolean hypercube, and each must contain at least one edge crossing the cut between $X$ and $Y$ (i.e. a monochromatic edge).
Hence $\D[f] \le n ( 2^{n-1} - 1 ) = \D[f^*_{n,n-1}]$, so $f$ is not \acm{}.
\end{proof}

\begin{lemma}\label{GLbound}
Let $n, d \in \N$, and let $g$ have maximal $\D[g]$ over all $(n-1,d)$-PTFs.
Then for every $(n,d)$-PTF $f$, $\displaystyle \D[f] \le \frac{2n}{n-1} \D[g]$.
\end{lemma}
\begin{proof}
Let $n, d \in \N$, and let $g$ be an $(n-1,d)$-PTF with $\D[g]$ maximal.
Let $f$ be an $(n,d)$-PTF.
Any restriction $f'$ of $f$ to a function on $n-1$ variables is also a degree-$d$ PTF, so $\D[f'] \le \D[g]$.
There are $2n$ such restrictions $f'$, and each dichromatic edge of $f$ appears in exactly $n-1$ of them.
Hence $(n-1) \D[f] \le 2n \D[g]$, from which the desired result follows immediately.
\end{proof}

\begin{lemma}\label{GLparity}
Let $n, d \in \N$ with $d < n$.
If $n$ and $d$ have the same parity, and $(n-1,d)$-\cms{} do not exist, then $(n,d)$-\cms{} do not exist.
\end{lemma}
\begin{proof}
Assume that no $(n-1,d)$-PTF is \acm{}.
If $n$ and $d$ have the same parity, then every restriction of $f^*_{n,d}$ is equivalent (with respect to $\sim$) to $f^*_{n-1,d}$.
There are $2n$ such restrictions, and each dichromatic edge of $f^*_{n,d}$ appears in exactly $n-1$ of them, so $\D[f^*_{n,d}] = \frac{2n}{n-1} \D[f^*_{n-1,d}]$.
Hence by \Lemming{} \ref{GLbound}, $(n,d)$-\cms{} do not exist.
\end{proof}

\begin{corollary}\label{GLd=n-2}
For every $n \in \N$, $(n,n-2)$-\cms{} do not exist.
\end{corollary}
\begin{proof}
This follows from Corollary \ref{GLd=n-1} and \Lemming{} \ref{GLparity}.
\end{proof}

\begin{corollary}\label{GLn<6}
Let $n, d \in \N$.
If $d \le n \le 5$ and $(n,d) \ne (5,2)$, then $(n,d)$-\cms{} do not exist.
\end{corollary}
\begin{proof}
This follows from Corollaries \ref{GLd=1}, \ref{GLd=n-1} and \ref{GLd=n-2}, and the fact that $(n,d)$-\cms{} trivially do not exist when $d \in \lbrace 0, n \rbrace$.
\end{proof}

\subsection{A Simple Counterexample}

In the statement of Corollary \ref{GLn<6}, the caveat $(n,d) \ne (5,2)$ cannot be removed.

\begin{lemma}
There exists a unique $(5,2)$-\cm{} $f_{5,2}$, modulo $\sim$.
\end{lemma}
\begin{proof}
In the case where $n=5$ and $d=2$, $p^*_{5,2}(x) = x(x - 2)$, and $\D[f^*_{5,2}] = 50$.
Let $q \in \R [ x, y ]$ be defined by $q ( x, y ) := 3y^2 - x^2 + 2xy + y - x - 3$, let $q' \in \R [ x_1, \ldots, x_5 ]$ such that $q' ( x_1, \ldots, x_5 ) := q ( x_1 + x_2, x_3 + x_4 + x_5 )$, and let $f_{5,2} := \sgn \circ q'$.
Since $q$ is quadratic, $f_{5,2}$ is a $(5,2)$-PTF.
It is not difficult to verify that $\D[f_{5,2}] = 51 > 50 = \D[f^*_{5,2}]$, so $f_{5,2}$ is a $(5,2)$-\cm{}.
For uniqueness, see Appendix \ref{appendix-5,2}.
\end{proof}

The existence of a $(5,2)$-\cm{} precludes the use of \Lemming{} \ref{GLparity} to prove that $(6,2)$-\cms{} do not exist.
However, the uniqueness of $f_{5,2}$, along with the fact that it only has one additional dichromatic edge, allows for a proof using \Lemming{} \ref{GLbound}.
\begin{lemma}\label{GL6,2}
For every $d$, $(6,d)$-\cms{} do not exist.
\end{lemma}
\begin{proof}
The cases $d \in \lbrace 0, 1, 4, 5, 6 \rbrace$ have already been covered.
For $d=3$, see Appendix \ref{appendix-6,3}.
The case $d=2$ remains.
Assume for the sake of contradiction that $f$ is a $(6,2)$-\cm{}.
The dichromatic count of every boolean function on an even number of variables is an even integer.
Since for every $(5,2)$-PTF $g$, $\D[g] \le 51$, \Lemming{} \ref{GLbound} implies that $120 < \D[f] \le 122.4$, and hence that $\D[f] = 122$.
There are $12$ restrictions of $f$ to a function $g$ on $5$ variables, all of which satisfy $\D[g] \le 51$.
Every dichromatic edge in $f$ appears in exactly five such $g$, so the expectation over a uniformly random restriction $g$ of $\D[g]$ is $\displaystyle \frac{5}{12} \cdot 122 > 50.5$.
Since $\D[g]$ is always an integer, $\D[g] = 51$ with probability strictly greater than $1/2$.
In particular, there exists $i$ such that $f|_{x_i = -1} \sim f|_{x_i = 1} \sim f_{5,2}$ (*).
However, it is easily verified (see Appendix \ref{appendix-6,2}) that no function $f$ satisfying both (*) and $\D[f] = 122$ is a $(6,2)$-PTF.
This contradicts the initial choice of $f$.
Hence no $(6,2)$-\cms{} exist.
\end{proof}

This also completes the proof of Theorem \ref{mainconfirm}. \hfill$\Box$

\subsection{Extension to Odd $n$}
We may extend $f_{5,2}$ to an $(n,2)$-\cm{} for any odd $n \ge 5$.

\begin{theorem}\label{extension}
For every odd $n \in \N$ with $n \ge 5$, there exists an $(n,2)$-\cm{} $f_{n,2}$ with\\$\displaystyle \D[f_{n,2}] \in \left( 1 + \Omega \left( n^{-1} \right) \right) \D[f^*_{n,2}]$.
\end{theorem}
Intuitively, $f_{n,2}$ behaves exactly as $f_{5,2}$, with the additional variables contributing to the second argument of $q$.\\
\begin{proof}
Let $n \ge 5$ be an odd integer.
Let $A := \{-2, 0, 2\}$ and $B := 2\Z + 1$.
Let $H$ be the $n$-dimensional boolean hypercube, and let $G$ be the graph with vertex set $A \times B$ and an edge between $u$ and $v$ exactly when $\lVert u-v \rVert_1 = 2$.
Let $\phi : H \rightarrow G$ be the graph homomorphism defined by $\phi ( x_1, \ldots, x_n ) := \left( x_1 + x_2, x_3 + \ldots + x_n \right)$.
Let $q ( x, y ) := 3y^2 - x^2 + 2xy + y - x - 3$ as above, let $f := \sgn \circ q$, and take $f_{n,2} := f \circ \phi$.
Note that because $\phi$ is a graph homomorphism, we may compute $\D[f_{n,2}]$ by counting the dichromatic edges $e$ induced by $f$ on $G$, weighted by $\phi^{-1}(e)$.
To this end, we observe that an edge $e$ between $(2i-2, 2j+2-n)$ and $(2i-2, 2j-n)$ has a preimage under $\phi$ of cardinality
\begin{align*}
\lvert \phi^{-1}(e) \rvert = {2 \choose i} {n-2 \choose j} \left( n - 2 - j \right).
\end{align*}
Similarly, for an edge $e$ between $(2i-2, 2j+2-n)$ and $(2i, 2j+2-n)$,
\begin{align*}
\lvert \phi^{-1}(e) \rvert = {2 \choose i} {n-2 \choose j} \left( 2 - i \right).
\end{align*}
We observe that $q$ is positive on $A \times B$ except at the four points $\{ (-2, 1), (0, -1), (2, -1), (2, 1) \}$.
Hence $f$ gives nine dichromatic edges, as indicated by the black lines below.
\begin{center}
\begin{tabular}{cr|cccccccc}
&\multicolumn{1}{r}{} & \multicolumn{8}{c}{$j$} \\
&& $2-n$ & $\cdots$ & $-3$ & $-1$ & $+1$ & $+3$ & $\cdots$ & $n-2$ \\
\cline{2-10}
& $-2$ & $\color{cyan} +$ & $\cdots$ & $\color{cyan} +$ & $\color{cyan} +$ & \multicolumn{1}{|c|}{$\color{red} -$} & $\color{cyan} +$ & $\cdots$ & $\color{cyan} +$ \\
\cline{6-7}
$i$ & $0$ & $\color{cyan} +$ & $\cdots$ & $\color{cyan} +$ & \multicolumn{1}{|c|}{$\color{red} -$} & $\color{cyan} +$ & $\color{cyan} +$ & $\cdots$ & $\color{cyan} +$ \\
\cline{7-7}
& $+2$ & $\color{cyan} +$ & $\cdots$ & $\color{cyan} +$ & \multicolumn{1}{|c}{$\color{red} -$} & \multicolumn{1}{c|}{$\color{red} -$} & $\color{cyan} +$ & $\cdots$ & $\color{cyan} +$
\end{tabular}
\end{center}
Summing the above expressions over these nine edges, we have
\begin{align*}
\D[f_{n,2}] &= {n-2 \choose \frac{n-1}{2}} \left( {2 \choose 0} n + {2 \choose 1} \left( n-1 \right) + {2 \choose 2} \left( n-1 \right) \right) \\
&= {n-2 \choose \frac{n-1}{2}} \left( 4n-3 \right) \\
&= {n-2 \choose \frac{n-1}{2}} \left( n-3 + 3n \right) \\
&= {n-2 \choose \frac{n+1}{2}} \left( n+1 \right) + {n-2 \choose \frac{n-1}{2}} 3n \\
&= \left( {n-1 \choose \frac{n+1}{2}} + {n-1 \choose \frac{n-1}{2}} \right) n + {n-2 \choose \frac{n+1}{2}} \\
&= {n \choose \frac{n+1}{2}} n + {n-2 \choose \frac{n+1}{2}} \\
&\in \left( 1 + \Theta(n^{-1}) \right) \D[f^*_{n,2}].
\end{align*}
Hence $f_{n,2}$ is an $(n,2)$-\cm{}, as desired.
\end{proof}

\subsection{The General Case}
Using a similar construction, we now prove the existence of \cms{} of arbitrary degree.
\begin{theorem}\label{general}
For every $n, d \in \N$ with $n \ge 7$ and $3 \le d \le n-3$, there exists an $(n, d)$-\cm{} $f_{n,d}$ with $\D[f_{n,d}] \in \left( 1 + \Omega \left( n^{-1} e^{-d^2/n} \right) \right) \D[f^*_{n,d}]$.
\end{theorem}
We first consider the case where $n$ and $d$ have the same parity.
The case where $n$ and $d$ have opposite parity is similar but handled later.
\begin{theorem}
For every $n, d \in \N$ with $3 \le d \le n-4$ and $n-d$ even, there exists an $(n, d)$-\cm{} $f_{n,d}$ with $\D[f_{n,d}] \in \left( 1 + \Omega \left( n^{-1} e^{-d^2/n} \right) \right) \D[f^*_{n,d}]$.
\end{theorem}
\begin{proof}
Let $n, d$ be integers of the same parity with $3 \le d \le n-4$.
Let $A := \{-3, -1, 1, 3\}$ and let $B := 2\Z + d + 1$.
Let $H$ be the $n$-dimensional boolean hypercube, and let $G$ be the graph with vertex set $A \times B$ and an edge between $u$ and $v$ exactly when $\lVert u-v \rVert_1 = 2$.
Let $\phi$ be the graph homomorphism defined by $\phi(x_1, \ldots x_n) := (x_1 + x_2 + x_3, x_4 + \ldots + x_n)$.
We now define four polynomials $p_1, p_2, p_3, p_4$ on $A \times B$ as follows:
\begin{align*}
p_1(x,y) &:= (y-1+d)(y+1-d) \\
p_2(x,y) &:= 1 - 2 \left( x (d-1) + y \right)^2 \\
p_3(x,y) &:= (y-3+d)(y-5+d) \cdots (y+5-d)(y+3-d) \\
p_4(x,y) &:= x(x+2)(x-2)(y-4+d)(y-6+d) \cdots (y+6-d)(y+4-d)
\end{align*}
Since $p_1 \in \Omega(p_2)$, there exists $\eps' > 0$ such that for every $v \in G$ with $p_1(v) \ne 0$, $|p_1(v)| > |2 \eps' p_2(v)|$.
Similarly, $p_3 \in \Omega(p_4)$, so there exists $\eps \in (0, \eps']$ such that for every $v \in G$ with $p_3(v) \ne 0$, $|p_3(v)| > |\eps p_4(v)|$.
For instance, we may take $\eps  = \eps' = \left( 4d \right)^{-d}$.
Take $p := \left( p_1 + \eps p_2 \right) \cdot p_3 - \eps^2 p_4$, take $g := \sgn \circ p$, and take $f_{n,d} := g \circ \phi$.
Since $p_1$ and $p_2$ have degree $2$, $p_3$ has degree $d-2$, and $p_4$ has degree $d$, $f_{n,d}$ is a polynomial threshold function of degree $d$.
Towards computing $\D[f_{n,d}]$, we first consider the relevant behaviors of $p_1, p_2, p_3, p_4$ separately.
All four are integer-valued (evaluations of) polynomials on the domain $A \times B$.
Both $p_2$ and $p_4$ are always odd, so in particular, are non-zero everywhere.
Firstly, $p_3$ is positive when $y > d-3$, is zero when $|y| \le d-3$, and has the same sign as $(-1)^d$ when $y < 3-d$.
Clearly, $p_1$ is positive when $|y| > d-1$ and zero when $|y| = d-1$.
By choice of $\eps$, $p_1 + \eps p_2$ is never in the interval $(-\eps, \eps)$ and always has the same sign as $p_1$ when $p_1$ is non-zero.
Similarly, $p$ is always non-zero and always has the same sign as $\left( p_1 + \eps p_2 \right) \cdot p_3$ when $p_3$ is non-zero.
Hence we may rewrite $g$ as the following piecewise function:
\begin{align*}
g(x,y) &= \left\lbrace
\begin{array}{ll}
(-1)^d & y < 1-d \\
\sgn ((-1)^d p_2(x,y)) & y = 1-d \\
\sgn (p_4(x,y)) & |y| < d-1 \\
\sgn (p_2(x,y)) & y = d-1 \\
1 & y > d-1
\end{array}
\right.
\end{align*}
Since $p_2$ is positive only at the two points $(-1, d-1)$ and $(1, 1-d)$ when $|y| = d-1$, the above piecewise representation shows that when $|y| \le d-1$, $g(x,y)$ computes the parity function except at the two points $(3, d-1)$ and $(-3, 1-d)$ (illustrated in Figures \ref{evenfigure} and \ref{oddfigure}).
\begin{figure}[h]
\begin{center}
\begin{tabular}{cr|ccccccccc}
& \multicolumn{1}{r}{} & \multicolumn{9}{c}{$y$} \\
&& $\cdots$ & $-1-d$ & $1-d$ & $3-d$ & $\cdots$ & $d-3$ & $d-1$ & $d+1$ & $\cdots$ \\
\cline{2-11}
\multirow{4}{*}{$x$} & $-3$ & $\cdots$ & $\color{cyan} +$ & \multicolumn{1}{|c}{$\color{red} -$} & $\color{red} -$ & \multicolumn{1}{|c|}{$\cdots$} & $\color{cyan} +$ & \multicolumn{1}{|c|}{$\color{red} -$} & $\color{cyan} +$ & $\cdots$ \\
\cline{6-9}
& $-1$ & $\cdots$ & $\color{cyan} +$ & \multicolumn{1}{|c|}{$\color{red} -$} & $\color{cyan} +$ & \multicolumn{1}{|c|}{$\cdots$} & $\color{red} -$ & \multicolumn{1}{|c}{$\color{cyan} +$} & $\color{cyan} +$ & $\cdots$ \\
\cline{5-9}
& $+1$ & $\cdots$ & $\color{cyan} +$ & \multicolumn{1}{c|}{$\color{cyan} +$} & $\color{red} -$ & \multicolumn{1}{|c|}{$\cdots$} & $\color{cyan} +$ & \multicolumn{1}{|c|}{$\color{red} -$} & $\color{cyan} +$ & $\cdots$ \\
\cline{5-8}
& $+3$ & $\cdots$ & $\color{cyan} +$ & \multicolumn{1}{|c|}{$\color{red} -$} & $\color{cyan} +$ & \multicolumn{1}{|c|}{$\cdots$} & $\color{red} -$ & \multicolumn{1}{c|}{$\color{red} -$} & $\color{cyan} +$ & $\cdots$
\end{tabular}
\end{center}
\caption{Illustration of $g$ in the case where $n$ and $d$ are both even}
\label{evenfigure}
\end{figure}
\begin{figure}[h]
\begin{center}
\begin{tabular}{cr|ccccccccc}
& \multicolumn{1}{r}{} & \multicolumn{9}{c}{$y$} \\
&& $\cdots$ & $-1-d$ & $1-d$ & $3-d$ & $\cdots$ & $d-3$ & $d-1$ & $d+1$ & $\cdots$ \\
\cline{2-11}
\multirow{4}{*}{$x$} & $-3$ & $\cdots$ & $\color{red} -$ & \multicolumn{1}{|c}{$\color{cyan} +$} & $\color{cyan} +$ & \multicolumn{1}{|c|}{$\cdots$} & $\color{cyan} +$ & \multicolumn{1}{|c|}{$\color{red} -$} & $\color{cyan} +$ & $\cdots$ \\
\cline{6-9}
& $-1$ & $\cdots$ & $\color{red} -$ & \multicolumn{1}{|c|}{$\color{cyan} +$} & $\color{red} -$ & \multicolumn{1}{|c|}{$\cdots$} & $\color{red} -$ & \multicolumn{1}{|c}{$\color{cyan} +$} & $\color{cyan} +$ & $\cdots$ \\
\cline{5-9}
& $+1$ & $\cdots$ & $\color{red} -$ & \multicolumn{1}{c|}{$\color{red} -$} & $\color{cyan} +$ & \multicolumn{1}{|c|}{$\cdots$} & $\color{cyan} +$ & \multicolumn{1}{|c|}{$\color{red} -$} & $\color{cyan} +$ & $\cdots$ \\
\cline{5-8}
& $+3$ & $\cdots$ & $\color{red} -$ & \multicolumn{1}{|c|}{$\color{cyan} +$} & $\color{red} -$ & \multicolumn{1}{|c|}{$\cdots$} & $\color{red} -$ & \multicolumn{1}{c|}{$\color{red} -$} & $\color{cyan} +$ & $\cdots$
\end{tabular}
\end{center}
\caption{Illustration of $g$ in the case where $n$ and $d$ are both odd}
\label{oddfigure}
\end{figure}
We now define $g' : A \times B \rightarrow \{ -1, 1 \}$ by $g'(\phi_{n,3}(x)) := -f^*_{n,d}(x)$.
Note that because $f^*_{n,d}$ is symmetric, this gives a well-defined function $g'$.
It is easily verified that for all $(x,y) \in A \times B$ such that $|y| \le d-1$ and at the two points $(3, -1-d)$ and $(-3, d+1)$, $g'(x,y) = g(x,y)$, and that for all other $(x,y) \in A \times B$, $g'(x,y) = -g(x,y)$.
Hence there are ten edges $\{ u, v \}$ in $G$ for which $g(u)g(v) \ne g'(u)g'(v)$.
This allows us to compute $\D[f_{n,d}]$ as follows:
\begin{align*}
\D[f_{n,d}] =& \ \D[f^*_{n,d}] \\
&+ (n+d-2) {n-3 \choose \frac{n-d-4}{2}} + 3(n+d-2) {n-3 \choose \frac{n-d-4}{2}} \\
&- 3(n+d-2) {n-3 \choose \frac{n-d-4}{2}} - 6 {n-3 \choose \frac{n-d-4}{2}} - (n-d-4) {n-3 \choose \frac{n-d-4}{2}} \\
=& \ \D[f^*_{n,d}] + (2d-4) {n-3 \choose \frac{n-d-4}{2}} \\
\in& \ \left( 1 + \Omega \left( \frac{{n \choose \frac{n-d}{2}}}{n {n \choose \frac{n}{2}}} \right) \right) \D[f^*_{n,d}] \\
\subseteq& \ \left( 1 + \Omega \left( n^{-1} e^{-d^2/n} \right) \right) \D[f^*_{n,d}]. &&\text{(Stirling's Inequality)}
\end{align*}
Hence $f_{n,d}$ is an $(n,d)$-\cm{}, as desired.
\end{proof}

\begin{theorem}
For every $n,d \in \N$ with $n \ge 7$, $3 \le d \le n-3$ and $n-d$ odd, there exists an $(n,d)$-\cm{} $f_{n,d}$ with $\displaystyle \D[f_{n,d}] \in \left( 1 + \Omega \left( n^{-1} e^{-d^2/n} \right) \right) \D[f^*_{n,d}]$.
\end{theorem}
\begin{proof}
The proof proceeds similarly to the previous case.
We define $A$, $B$, $H$, $G$, $p_1$, $p_2$, $p_3$, $p_4$, $p$, $g$, and $g'$ as above, and we define $\psi ( x_1, \ldots, x_n ) := x_1 + x_2 + x_3, 1 + x_4 + \ldots + x_n )$.
We now define $f_{n,d} := g \circ \psi$ analogously to above.
The computation of $\D[f_{n,d}]$ now proceeds as follows:
\begin{align*}
\D[f_{n,d}] =& \ \D[f^*_{n,d}] \\
&+ \frac{n+d-3}{2} {n-3 \choose \frac{n-d-3}{2}} + 3 \frac{n+d-3}{2} {n-3 \choose \frac{n-d-3}{2}} \\
&- 3 \frac{n+d-3}{2} {n-3 \choose \frac{n-d-3}{2}} - 3 {n-3 \choose \frac{n-d-3}{2}} - \frac{n-d-3}{2} {n-3 \choose \frac{n-d-3}{2}} \\
&+ \frac{n+d-1}{2} {n-3 \choose \frac{n-d-5}{2}} + 3 \frac{n+d-1}{2} {n-3 \choose \frac{n-d-5}{2}} \\
&- 3 \frac{n+d-1}{2} {n-3 \choose \frac{n-d-5}{2}} - 3 {n-3 \choose \frac{n-d-5}{2}} - \frac{n-d-5}{2} {n-3 \choose \frac{n-d-5}{2}} \\
=& \ \D[f^*_{n,d}] + (d-3) {n-3 \choose \frac{n-d-3}{2}} + (d-1) {n-3 \choose \frac{n-d-5}{2}} \\
\in& \ \left( 1 + \Omega \left( \frac{{n \choose \frac{n-d}{2}}}{n {n \choose \frac{n}{2}}} \right) \right) \D[f^*_{n,d}] \\
\subseteq& \ \left( 1 + \Omega \left( n^{-1} e^{-d^2/n} \right) \right) \D[f^*_{n,d}]. &&\text{(Stirling's Inequality)}
\end{align*}
Hence $f_{n,d}$ is an $(n,d)$-\cm{}, as desired.
\end{proof}

This also completes the proofs of Theorems \ref{general} and \ref{mainrefute}. \hfill$\Box$

\section{Conclusion}

For almost all $d$ and almost all $n$, we refute the Gotsman-Linial Conjecture (Conjecture \ref{GL}) with a multiplicative separation of $1 + \Theta_d \left( n^{-1} \right)$.
This separation is too weak to refute most known applications of the conjecture.
We would need to improve $1 + \Theta_d \left( n^{-1} \right)$ to $\omega(1)$ to refute the Asymptotic Gotsman-Linial Conjecture (Conjecture \ref{GLweak}), on which the applications depend.
Although for every $(n,d)$-\cm{} $f$ given in this paper, $\D[f] > \D[f^*_{n,d}]$, it should be noted that the RHS is still an upper bound in a \emph{limiting} sense.
This, along with the intuition presented in Section \ref{intuition}, invites the following revised conjecture.
\begin{conjecture}[Gotsman-Linial - Limit]\label{GLlimit}
Let $f : \booleanfn$ be an $(n,d)$-PTF.
Then the average sensitivity $\AS[f] \le d \AS[f^*_{n,1}]$.
\end{conjecture}

Conjecture \ref{GLlimit} would resolve the remaining cases of Conjectures \ref{GLcm} and \ref{GL}, i.e.
\begin{conjecture}
For every even $n$, $(n,2)$-\cms{} do not exist.
\end{conjecture}

Furthermore, our revised conjecture would imply the Asymptotic Gotsman-Linial Conjecture (Conjecture \ref{GLweak}) and its consequent applications.

\section*{Acknowledgments}
The author would like to thank Ryan Williams especially for inspiration, advice, feedback, and an admirable tolerance of cheesemonkeys; the Williams family, the Chap-people, Henry Qin, and Carolyn Kim for moral support and a good work environment; and Not Luke the goldfish for surviving.

\newcommand{\etalchar}[1]{$^{#1}$}

\appendix

\section{Appendix}\label{appendix-computer}

Here we describe how a computer search resolved the cases of $n=6$ and $d=2,3$ of the Gotsman-Linial Conjecture.
First, we note that the problem of determining whether a boolean function $f : \booleanfn$ is an $(n,d)$-PTF is equivalent to determining whether a particular linear program has any feasible solution.
Unfortunately, leveraging this fact to compute the maximal average sensitivity of an $(n,d)$-PTF with a na\"ive exhaustive search takes doubly exponential time so is intractable for large $n$ (i.e. $n > 4$).
However, by using \Lemming{} \ref{GLbound}, we can conduct a more efficient search.
We maintain a partial function $f$ and conduct a DFS in which we define $f$ successively on inputs in increasing order of Hamming weight.
This allows us to keep bounds on $\D[f]$ by counting the edges that are already constrained to be monochromatic or dichromatic.
When $\D[f]$ becomes too low or too high, we can prune the search and backtrack before fully defining $f$, allowing (tolerably) efficient searches up to $n=6$.

\subsection{$(n,d) = (5,2)$}\label{appendix-5,2}
In the case $(n,d) = (5,2)$, \Lemming{} \ref{GLbound} implies that for every $(5,2)$-\cm{} $f$, we have $50 < \D[f] \le 60$.
Because a random boolean function on $5$ variables has far fewer than $50$ dichromatic edges with high probability, most search branches are pruned early.
The modified search confirmed that every $(5,2)$-\cm{} $f$ with $50 < \D[f] \le 60$ satisfies $f \sim f_{5,2}$, and hence that $f_{5,2}$ is the unique $(5,2)$-\cm{}.

\subsection{$(n,d) = (6,2)$}\label{appendix-6,2}
The proof of \Lemming{} \ref{GL6,2} relies on the claim that for every $(6,2)$-PTF $f$ and for every variable $x_i$, either $f|_{x_i = -1} \not \sim f_{5,2}$, $f|_{x_i = 1} \not \sim f_{5,2}$, or $\D[f] \ne 122$.
By symmetry, we may remove the dependence on $i$, and we may turn one of the equivalences into an equality.
The above is equivalent to the claim that for every $(6,2)$-PTF $f$, either $f|_{x_1 = -1} \not \sim f_{5,2}$, $f|_{x_1 = 1} \ne f_{5,2}$, or $\D[f] \ne 122$.
Because a restriction of a $(6,2)$-PTF is a $(5,2)$-PTF and must therefore have no more than $51$ dichromatic edges, it suffices to show that for every function $f : \lbrace -1, 1 \rbrace^6 \rightarrow \lbrace -1, 1 \rbrace$ such that $f|_{x_1 = -1} \sim f|_{x_1 = 1} = f_{5,2}$ and $\D[f] = 122$, there exist $i$ and $b$ such that $\D[f|_{x_i = b}] > 51$.
There are fewer than $1000$ functions $f$ satisfying both $f|_{x_1 = -1} \sim f|_{x_1 = 1} = f_{5,2}$ and $\D[f] = 122$ (and these are easily enumerated), and each has only $10$ relevant restrictions to five variables, so the last claim is easily verified with a quick computer search.

\subsection{$(n,d) = (6,3)$}\label{appendix-6,3}
In the case $(n,d) = (6,3)$, \Lemming{} \ref{GLbound} implies that for every $(6,3)$-\cm{} $f$, we have $150 < \D[f] \le 168$.
Because a random boolean function on $6$ variables has far fewer than $150$ dichromatic edges with overwhelming probability, most search branches are pruned early.
The modified search confirmed that no $(6,3)$-PTF $f$ satisfying $150 < \D[f] \le 168$ is \acm{}, and hence that no $(6,3)$-\cms{} exist.


\begin{thebibliography}{BCG{\etalchar{+}}96}

\bibitem[ABFR94]{ABFR}
James Aspnes, Richard Biegel, Merrick Furst, and Steven Rudich.
\newblock The Expressive Power of Voting Polynomials.
\newblock \textit{Combinatorica}, 14(2):1-14, 1994.

\bibitem[Bou79]{Bourgain}
Jean Bourgain.
\newblock Walsh subspaces of $L^p$ product spaces.
\newblock \textit{S\'eminaire D'Analyse Fonctionelle}, \'Ecole Polytechnique, Centre de Mathematiques, pp. IV.1-IV.9, 1979.

\bibitem[CSS16]{CSS}
Ruiwen Chen, Rahul Santhanam, and Srikanth Srinivasan.
\newblock Average-Case Lower Bounds and Satisfiability Algorithms for Small Threshold Circuits.
\newblock \textit{Proceedings for the $31$st Conference on Computational Complexity}, pp. 1:1-1:35, 2016.

\bibitem[DHKMRST10]{souplite}
Ilias Diakonikolas, Prahladh Harsha, Adam Klivans, Raghu Meka, Prasad Raghavendra, Rocco A. Servedio, and Li-Yang Tan.
\newblock Bounding the average sensitivity and noise sensitivity of polynomial threshold functions.
\newblock In \textit{ACM Symposium on Theory of Computing (STOC)}, pp. 533--542, 2010.

\bibitem[DRST14]{DRST}
Ilias Diakonikolas, Prasad Raghavendra, Rocco A. Servedio, and Li-Yang Tan.
\newblock Average sensitivity and noise sensitivity of polynomial threshold functions.
\newblock SIAM Journal on Computing, pp. 231--253, 2014.

\bibitem[FHHMOSWW14]{alphabetsoup}
Yuval Filmus, Hamed Hatami, Steven Heilman, Elchanan Mossel, Ryan O’Donnell, Sushant Sachdeva, Andrew Wan, and Karl Wimmer.
\newblock Real Analysis in Computer Science: A collection of Open Problems.
\newblock Simons Institute, Berkeley, CA, compiled in 2014. 
\newblock URL: https://simons.berkeley.edu/sites/default/files/openprobsmerged.pdf

\bibitem[GS10]{GS}
Parikshit Gopalan and Rocco Servedio.
\newblock Learning and Lower Bounds for $\AC^0$ with Threshold Gates.
\newblock \textit{Proceedings for the $14$th International Workshop on Randomization and Computation}, pp. 588-601, 2010.

\bibitem[GL94]{GL94}
Craig Gotsman and Nathan Linial.
\newblock Spectral Properties of Threshold Functions.
\newblock \textit{Combinatorica}, 14(1):35-50, 1994.

\bibitem[HKM09]{HKM}
Prahladh Harsha, Adam Klivans and Raghu Meka. 
\newblock Bounding the sensitivity of polynomial threshold functions.
\newblock \textit{arXiv:0909.5175}, 2009.

\bibitem[Kan12]{K12}
Daniel M. Kane.
\newblock A structure theorem for poorly anticoncentrated Gaussian chaoses and applications to the study of polynomial threshold functions.
\newblock In \textit{Foundations of Computer Science (FOCS)}, pp. 91--100, 2012.

\bibitem[Kan13]{K13}
Daniel M. Kane.
\newblock The Correct Exponent for the Gotsman-Linial Conjecture.
\newblock \textit{arXiv:1210.1283}, 2013.

\bibitem[KW16]{KW}
Daniel Kane and Ryan Williams.
\newblock Super-Linear Gate and Super-Quadratic Wire Lower Bounds for Depth-Two and Depth-Three Threshold Circuits.
\newblock \textit{Proceedings for the $48$th Annual ACM SIGACT Symposium on the Theory of Computing}, pp. 633-643, 2016.

\bibitem[LMN93]{LMN}
Nathan Linial, Yishay Mansour, and Noam Nisan.
\newblock Constant depth circuits, Fourier transform, and learnability.
\newblock \textit{Journal of the ACM}, 40(3):607-620, 1993.

\bibitem[MP68]{MP}
Marvin Minsky and Seymour Papert.
\newblock \textit{Perceptrons: an Introduction to Computational Geometry (Expanded Edition)}.
\newblock MIT Press, Cambridge, MA, 1988.

\bibitem[OD12]{OD12}
Ryan O'Donnell.
\newblock Open problems in analysis of Boolean functions.
\newblock \textit{arXiv preprint arXiv:1204.6447}, 2012.

\bibitem[OD14]{OD14}
Ryan O'Donnell.
\newblock \textit{Analysis of Boolean Functions}.
\newblock Cambridge University Press, New York, New York, 2014.

\bibitem[ON71]{ON}
Patrick E. O'Neil.
\newblock Hyperplane Cuts of an $n$-Cube.
\newblock \textit{Discrete Maths}. \textbf{1} (1971), 193-195.
\end{thebibliography}
\end{document}